\documentclass[12pt]{article}
\usepackage[reqno]{amsmath}
\usepackage{amssymb, theorem, enumerate}
\usepackage{graphicx}

\theoremstyle{change}
{\theorembodyfont{\slshape}
\newtheorem{theorem}{Theorem.}[section]
\newtheorem{lemma}[theorem]{Lemma.}
}
\theorembodyfont{\rmfamily}

\newcommand\lref[1]{Lemma~\ref{lem:#1}}
\newcommand\tref[1]{Theorem~\ref{thm:#1}}
\newcommand\cref[1]{Corollary~\ref{cor:#1}}

\def\proof{\noindent{{\sl Proof. }}}
\def\sqr#1#2{{\vbox{\hrule height.#2pt
    \hbox{\vrule width.#2pt height#1pt \kern#1pt
        \vrule width.#2pt}\hrule height.#2pt}}}
\def\eqed{\sqr53}
\def\qed{%
    \ifmmode\eqno\eqed
    \else\nobreak\ \hfill\eqed\medbreak\fi}

\newcommand\de{\delta}

\newcommand\cprod{\mathbin{\scriptscriptstyle\square}}
\newcommand\aut[1]{{\rm Aut}(#1)}

\DeclareMathOperator{\Alt}{Alt}
\newcommand\alt[1]{\Alt(#1)}

\title{Distinguishing Primitive Permutation Groups}
\author{
	Chris Godsil\\
	Combinatorics and Optimization\\
	University of Waterloo\\[3pt]
	\texttt{cgodsil@uwaterloo.ca}}

\begin{document}
\maketitle

\begin{abstract}
Let $G$ be a permutation group acting on a set $V$.  A partition $\pi$ of $V$
is \textsl{distinguishing} if the only element of $G$ that fixes each cell of $\pi$
is the identity.  The \textsl{distinguishing number} of $G$ is the minimum number
of cells in a distinguishing partition.  We prove that if $G$ is a primitive permutation 
group and $|V|\ge336$, its distinguishing number is two.
\end{abstract}

\section{Introduction}

Let $G$ be a permutation group acting on a set $V$.  A partition $\pi$ of $V$
is \textsl{distinguishing} if the only element of $G$ that fixes each cell of $\pi$
is the identity.  The \textsl{distinguishing number} of $G$ is the minimum number
of cells in a distinguishing partition.  The study of this parameter was introduced by
Albertson and Collins, see \cite{alco}.
For other recent work on this problem see, for example, the work of M. Chan \cite{chan-max}.

The result of this paper is that if $G$ is primitive and on $V$ and $|V|\ge336$, 
then its distinguishing number is two.

We note that Theorem 1 from Seress \cite{seress} yields the same conclusion with
an upper bound of $32$.  The best we can say for this paper is that our methods are 
more elementary, and may have some independent interest.

\section{A Lower Bound}

Let $G$ be a permutation group.  The \textsl{degree} of $G$ is just the size of $V$.
The \textsl{minimum degree} of $G$ is the minimum number of points moved by
a non-identity element of $G$.

We need the following lower bound on the order of a permutation group in terms of its
minimum degree.  The argument is based on Cameron \cite{pjc-power}.

\begin{lemma}
\label{lem:pjc}
If $G$ is a permutation group with distinguishing number at least three and the
minimum degree of $G$ is $\de$, then $|G| \ge1+2^{\de/2}$.
\end{lemma}

\proof
Assume $v=|V|$.
We consider the action of $G$ as a group of permutations of the power set of $V$.
We recall that the number of orbits of $G$ on subsets of $V$ is equal to the average number
of subsets fixed by an element of $G$.  If $g\in G$ and $g$ has exactly $k$ orbits,
then it fixes exactly $2^{k}$ subsets.  Since the minimum degree of $G$ is $\de$, we see 
that $g$ has at most $v-\de$ orbits of length one, and at most $\de/2$ orbits of length 
greater than one.  Therefore $g$ fixes at most $2^{v-\frac{\de}2}$ subsets of $V$, 
and the number of orbits of $G$ on subsets is at most
\[
\frac1{|G|}(2^{v}+(|G|-1)2^{v-\frac{\de}2}).
\]

If the distinguishing number of $G$ is greater than two, then any non-empty proper
subset of $V$ must be fixed by at least one non-identity element of $G$, and so the
orbit of any subset has size at most $|G|/2$.  Consequently the number of orbits of
$G$ is at least
\[
\frac{2^{v}}{|G|/2} =\frac{2^{v+1}}{|G|}.
\]
Combining these two bounds and rearranging yields the lemma.\qed

\section{Primitive Groups}

The \textsl{Johnson graph} $J(m,\ell)$ has the $\ell$-subsets of $\{1,\ldots,m\}$
as its vertices, and two $\ell$-subsets are adjacent if their intersection has size $\ell-1$.
The graphs $J(m,\ell)$ and $J(m,m-\ell)$ are isomorphic, and so we will
assume that $2\ell\le m$.  The complete graphs are Johnson graphs (with $\ell=1$).

A subset $S$ of $G$ is a \textsl{block of imprimitivity} for $G$ if for each element 
$g$ of $G$,
either $S^{g}=S$ or $S\cap S^{g}=\emptyset$.  Trivial examples of blocks are $\emptyset$, 
any singleton from $V$ or $V$ itself.  If $G$ is transitive and there are no other blocks 
we say 
that $G$ is \textsl{primitive}, otherwise it is \textsl{imprimitive}.  The size of a 
non-empty block must divide $|V|$, whence we see that if $G$ is transitive and its degree 
is a prime number, then $G$ is primitive.

The \textsl{Cartesian product} $X\cprod Y$ of graphs $X$ and $Y$ has vertex
set $V(X)\times V(Y)$, and vertices $(x_1,y_1)$, $(x_2,y_2)$ in the product
are adjacent if $x_1=y_1$ and $x_2$ is adjacent to $y_2$, or if
$x_2=y_2$ and $x_1$ is adjacent to $y_1$.  We use $X^{\cprod\,n}$
to denote the $n$-th Cartesian power of $X$, the Cartesian product of $n$ copies
of $X$.

For our purposes the main result (Theorem~1) of Guralnick and Magaard \cite{guma} can be summarized as follows.

\begin{theorem}
\label{thm:gurmag}
Let $G$ be a primitive permutation group of degree $v$.  If the minimal degree of $G$
is at most $v/2$, then one the following holds:
\begin{enumerate}[(a)]
\item
$G$ is affine over $GF(2)$ and its minimal degree is $v/2$.
\item
$G$ is a transitive subgroup of the automorphism group of the Cartesian power
of a Johnson graph.
\end{enumerate}
\end{theorem}

The following result comes from Mar\'oti \cite[Theorem~1.1]{maroti}.

\begin{theorem}
\label{thm:maroti}
Let $G$ be a primitive permutation group of degree $v$.  Then one of the following holds:
\begin{enumerate}[(a)]
\item
$G$ is a transitive subgroup of the automorphism group of the Cartesian power of a
Johnson graph.
\item
$v\in\{11,12,23,24\}$ and $G$ is one of the Mathieu groups.
\item
$|G| < v^{1+\lceil\log_{2}v\rceil}$.\qed
\end{enumerate}
\end{theorem}

\section{Distinguishing Primitive Groups}
\label{sec:DPGs}

Imrich and Klavzar prove the following result in \cite[Theorem~1.1]{imkl}.

\begin{theorem}
	\label{thm:imkl}
The distinguishing  number of a Cartesian power of connected graph is two,
except in the cases $K_{2}^{\cprod\,2}$, $K_{3}^{\cprod\,2}$ and $K_{2}^{\cprod\,3}$
(when it is three).
\end{theorem}

If the graph $X$ has $v$ vertices then 
\[
\aut{X^{\cprod\,n}} \le\aut{K_{v}^{\cprod\,n}}
\]
and therefore the distinguishing number of $X^{\cprod\,n}$ is bounded by that of 
$K_{v}^{\cprod\,n}$.

\begin{lemma}
	\label{lem:jvk}
If $m\ge2\ell$ then the distinguishing number of $J(m,\ell)$ is two unless $\ell=1$
or $\ell=2$ and $m=5$.
\end{lemma}

\proof
Suppose $\ell=2$.  Let $T_{m}$ be the graph constructed from the path on $n$ vertices
by adding an edge joining the second and fourth vertices.  If $m>5$, then $T_{m}$
is asymmetric and its line graph is asymmetric.  It follows that the partition of 
the vertices of $L(K_{m})$ into the vertices which represent edges of $T_{m}$
and the vertices which do not is a distinguishing partition for $\aut{J(m,2)}$ when $m\ge6$.

Now suppose $\ell\ge3$.  Let $HP_{m}$ denote the hypergraph with vertex set
$\{1,\ldots,2m\}$ and edges consisting of the $m+1$ triples of consecutive integers
from the vertex set.  Add the edge consisting of the integers $\{2,\ldots,m,m+2\}$.
This modified hyperpath is asymmetric, and it follows that $\aut{J(m,\ell)}$ has
distinguishing number two.\qed

\begin{theorem}
Suppose $G$ is a primitive permutation group of degree $v$ that does not contain $\alt{v}$.
If $v\ge 336$, then the distinguishing number of $G$ is two.
\end{theorem}

\proof
Suppose that the minimum degree of $G$ is at least $v/2$.  Then combining \lref{pjc}
and \tref{maroti} we have
\[
2^{v/4} <v^{1+\lceil\log_{2}v\rceil}.
\]
Taking logarithms of each side of this yields
\[
v/4 <(1+\lceil\log_{2}v\rceil) \log_{2}v.
\]
This inequality fails if $v\ge 336$, and therefore the theorem holds when
the minimum degree of $X$ is at least $v/2$.

If the minimum degree of $G$ is less than $v/2$, then by \tref{gurmag}
$G$ must be a subgroup of the automorphism group of the Cartesian power of a 
Johnson graph.  By \tref{imkl}, the distinguishing number of the $n$-th Cartesian power
of $J(m,\ell)$ is two unless $J(m,\ell)$ is $K_{2}$ or $K_{3}$ and $n\le3$,
or unless $n=1$.  Using \lref{jvk} we see that in all these cases $v<336$.
Since the Mathieu groups have degree at most 24, these can also be ignored.\qed

\section{Comments}

Both the theorem of Guralnick and Magaard and the theorem of Ma\'roti depend on the
classification of finite simple groups.  We do have the following important bound due to 
Babai which can be used in place of the theorem of Guralnick and Magaard, and whose 
proof is entirely combinatorial. (For details, see \cite[Theorem 5.3B]{dimo}.)

\begin{theorem}
\label{thm:babai}
If $G$ is a primitive permutation group of degree $v$ that does not contain $\alt{v}$,
then its minimum degree is greater than $\sqrt{v}/2$.
\end{theorem}

However using this yields a lower bound of order $2^{20}$, much larger than the bound of
336 that we computed.


\end{document}